\documentclass{agtart_a}
\pdfoutput=1
\usepackage[all]{xy}

\title[Combinatorial groups, Hopf invariant and the
exponent problem]{Applications of combinatorial groups to\\Hopf invariant and the
exponent problem} 

\author{Jelena Grbic}
\givenname{Jelena}
\surname{Grbi\'c}
\address{Department of Mathematical Sciences\\
   University of Aberdeen\\\newline
   Meston Building\\
   Aberdeen AB24 3UE\\UK}
\email{jelena@maths.abdn.ac.uk}
\urladdr{}

\author{Jie Wu}
\givenname{Jie}
\surname{Wu}
\address{Department of Mathematics\\
  National University of Singapore\\\newline
  2 Science Drive 2\\
  Singapore}
\email{matwuj@nus.edu.sg}
\urladdr{}

\volumenumber{6}
\issuenumber{}
\publicationyear{2006}
\papernumber{78}
\startpage{2229}
\endpage{2256}

\doi{}
\MR{}
\Zbl{}


\received{23 October 2006}
\revised{}
\accepted{24 October 2006}
\published{29 November 2006}
\publishedonline{29 November 2006}
\proposed{}
\seconded{}
\corresponding{}
\editor{CPR}
\version{}

\arxivreference{math.AT/0602204}

\let\xysavmatrix\xymatrix
\def\xymatrix{\disablesubscriptcorrection\xysavmatrix}
\AtBeginDocument{\let\bar\wbar\let\tilde\wtilde\let\hat\what}
\numberwithin{equation}{section}
\def\strutt{\vrule width 0pt depth 5pt height 12pt}

\def\co{\mskip .7mu\colon\thinspace}

\makeatletter
\def\cnewtheorem#1[#2]#3{\newtheorem{#1}{#3}[section]
\expandafter\let\csname c@#1\endcsname\c@thm}
\makeatother

\newtheorem{thm}{Theorem}[section] 
\cnewtheorem{lem}[thm]{Lemma}       
\cnewtheorem{proposition}[thm]{Proposition}
\cnewtheorem{corollary}[thm]{Corollary}

\theoremstyle{definition}
\cnewtheorem{defn}[thm]{Definition} 
\makeautorefname{defn}{Definition}
\newtheorem*{rem}{Remark}  
\newtheorem*{conj}{The Barratt Conjecture}
\newtheorem*{conjs}{The Strong Barratt Conjecture}

\newcommand{\lra}{\longrightarrow}
\newcommand{\Ker}{\mathop{\mathrm{Ker}}}
\renewcommand{\Im}{\mathop{\mathrm{Im}}}

\newcommand{\Hom}{\ensuremath{\mathrm{Hom}}}
\newcommand{\Coalg}{\ensuremath{\mathrm{Coalg}}}

\def\Id{\mathop{\rm Id}\nolimits}

\renewcommand{\H}{\ensuremath{\mathcal{H}}}

\def\Lie{\mathop{\rm Lie}\nolimits}
\def\tr{\mathop{\rm tr}\nolimits}
\newcommand{\dr}[3]{\ensuremath{#1\stackrel{#2}
{\longrightarrow}#3}}
\newcommand{\ddr}[5]{\ensuremath{#1\stackrel{#2}
{\longrightarrow}#3\stackrel{#4}{\longrightarrow}#5}}
\newcommand{\dddr}[7]{\ensuremath{#1\stackrel{#2}
{\longrightarrow}#3\stackrel{#4}{\longrightarrow}#5
\stackrel{#6}{\longrightarrow}#7}}
\newcommand{\ddddr}[9]{\ensuremath{#1\stackrel{#2}
{\longrightarrow}#3\stackrel{#4}{\longrightarrow}#5
\stackrel{#6}{\longrightarrow}#7}\stackrel{#8}{\longrightarrow}#9}

\def\Lie{\mathop{\rm Lie}\nolimits}
\def\id{\mathop{\rm id}\nolimits}

\def\hocolim{\mathop{\rm hocolim}\nolimits}

\def\pinch{\mathop{\rm pinch}\nolimits}

\begin{document}

\begin{htmlabstract}
<p class="noindent">
Combinatorial groups together with the groups of natural coalgebra
transformations of tensor algebras are linked to the groups of
homotopy classes of maps from the James construction to a loop
space. This connection gives rise to applications to homotopy
theory. The Hopf invariants of the Whitehead products are studied and
a rate of exponent growth for the strong version of the Barratt
Conjecture is given.
</p>
<h5>Corrections</h5>
<p class="noindent">
12 May 2007:  Minor corrections have been made to the statements of
Theorems 1.2 and 3.5 and to the paragraph above Lemma 4.2.
</p>
\end{htmlabstract}

\begin{abstract} 
Combinatorial groups together with the groups of natural coalgebra
transformations of tensor algebras are linked to the groups of
homotopy classes of maps from the James construction to a loop
space. This connection gives rise to applications to homotopy
theory. The Hopf invariants of the Whitehead products are studied and
a rate of exponent growth for the strong version of the Barratt
Conjecture is given.
\end{abstract}

\maketitle

\section{Introduction}\label{sec1}
The study of the groups of homotopy classes of maps from one
topological space to another has always been the central problem of
algebraic topology. In this paper we are concerned with the natural
maps from loop suspensions to loop spaces. To study them we develop
a method which arises from the relations between combinatorial
groups and the natural coalgebra transformations of tensor algebras
established in the predecessor~\cite{GW} to this paper. In
particular, let $R$ be a commutative ring with identity, and let
$\Coalg(A(-),B(-))$ denote the group of natural coalgebra
transformations over $R$ between functors $A$ and $B$. The tensor
algebra $T(V)$ generated by a free $R$--module $V$ has a natural
coalgebra filtration (the James filtration) $\{J_n(V)\}_{n\geq 0}$
given by $J_n(V)=\bigoplus_{j\leq n}T_j(V)$ for $n\geq 0$, where
$T_j(V)$ denotes the $j$-th stage of the tensor length filtration
for $T(V)$, that is, $T_j(V)=V^{\otimes j}$. Let $C(V)=J_1(V)$. For
two $R$--modules $C$ and $D$, define their smash product $C\wedge D$
to be the quotient module
\[
C\wedge D=(C\otimes D)/ (C\otimes_RR\oplus R\otimes_RD).
\]
By $V^{\wedge n}$ we denote the $n$--fold self smash product of $V$.

In~\cite{GW} we defined combinatorial groups $K^R_n$, $\H^R_n$,
$\negthinspace^R\H^{(l)}_n$, $\negthinspace^R\H^{(l),(k)}_n$ and
$K^R_n(k)$ such that there are group isomorphisms:
\[
\begin{array}{ll}
\Coalg(C(-)^{\otimes n}),T(-))\cong K_n^R &
\Coalg(J_n(-),T(-))\cong\H^R_n\\
\Coalg(J_n(-^{\otimes l}),T(-))\cong\,^R\H^{(l)}_n &
\Coalg(J_n(-^{\otimes l}),T(-^{\otimes
k}))\cong\,^R\H^{(l),(k)}_n\\
\Coalg(C(-)^{\otimes n}),T(-^{\wedge k}))\cong K_n^R(k) &
\end{array}
\]
for $1\leq n\leq\infty$. As a consequence of these group
identifications and the fact that $\Coalg(A(-),B(-))$ is contained
inside the algebra of natural linear transformations
$\Hom_R(A(-),B(-))$, there exist faithful representations of the
above combinatorial groups in the algebras of natural linear
transformations of tensor algebras.
\medskip

The first issue we address in this paper is a geometrical
realization of the stated algebraic relations. Restrict the
commutative ring $R$ to $\Z$ or $\Z/p^r$. Let $X$ be a space such
that its reduced diagonal $\dr{\bar{\Delta}\co X}{}{X\wedge X}$
is null homotopic. Denote the $n$--fold self smash product of $X$ by
$X^{(n)}$. The connection between the combinatorial groups $K^R_n$,
$\H^R_n$, $^R\H_n^{(l)}$, $^R\H_n^{(l)(k)}$ and $K^R_n(k)$ and the
groups of homotopy classes of maps from the James construction on
spaces with a null homotopic reduced diagonal to a loop space will
be given by constructing injective group homomorphisms:
\[
\begin{array}{ll}
e_X\co K^R_n\lra[X^{n},J(X)] &\quad e_X\co \H^R_n\lra\,
[J_n(X),J(X)]\\
e_X\co^R\H^{(l)}_n\lra[J_n(X^{(l)}),J(X)] &
\quad e_X\co\,^R\H^{(l)(k)}_n\lra[J_n(X^{(l)}),J(X^{(k)})]\\
e_X\co K^R_n(k)\lra[X^{n},J(X^{(k)})] &
\end{array}
\]
(see
Propositions~\ref{prop:georealK_n},~\ref{prop:georealK_n(k)},~\ref{prop:georealH_n},~\ref{prop:georealH_n^l},~\ref{prop:georealH_n^lk}).
\medskip

We then proceed to apply the just established group-geometrical
model to problems in homotopy theory. In 1931 Hopf defined what is
nowadays known as the Hopf invariant in order to study maps between
spheres of different dimensions which cannot be distinguished
homologically. Ever since, Hopf invariants and their relations with
Whitehead products have been widely studied (see for
instance Cohen and Taylor~\cite{CT} and Wu~\cite{Wu}), and these have various applications
in homotopy theory to the homotopy groups of spheres, exponent
problems and LS--category. The first application of our
group-geometrical model is to calculate the generalized $k$-th Hopf
invariant of the $n$-fold Whitehead product when $k$ does not divide
$n$ and to determine the second Hopf invariant of the $4$--fold
Whitehead product.
\medskip

Recall that for any pointed space $X$, the James--Hopf map $H_k\co
J(X)\lra J(X^{(k)})$ is combinatorially defined by
\[
H_k(x_1x_2\ldots x_n)=\prod_{1\leq i_1<i_2<\cdots <i_k\leq
n}(x_{i_1}x_{i_2}\ldots x_{i_k})
\]
with right lexicographical order in the product. The $n$--fold
Samelson product $\widetilde{W}_n$ on $X$ is given by the composite
\[
\ddr{\widetilde{W}_n\co X\wedge\ldots\wedge
X}{E\wedge\ldots\wedge E}{\Omega\Sigma X\wedge\ldots\wedge
\Omega\Sigma X}{[[\, ,\, ],\ldots,]}{\Omega\Sigma X}
\]
where $E\co X\lra\Omega\Sigma X$ is the canonical inclusion and
the second map $[[\, ,\, ],\ldots,]$ is the $n$--fold commutator. The
$n$--fold Whitehead product $W_n$ on $X$ is defined as the adjoint of
the $n$--fold Samelson product $\widetilde{W}_n$.

\begin{thm}
\label{HWknotn} Let $X$ be a pointed space with the null homotopic
reduced diagonal $\bar{\Delta}\co X\lra X\wedge X$. Then for
$n>k$,
\[
\ddr{J(X^{(n)})}{\Omega W_n}{J(X)}{H_k}{J(X^{(k)})}
\]
is null homotopic if $k$ does not divide $n$.
\end{thm}
The case when $k$ divides $n$ is much more subtle. Here we determine
the first non-trivial case, that is, the second Hopf invariant of
the $4$--fold Whitehead product, leaving the description of the
general case as a future project.
\begin{thm}
\label{HWkn}For $X$ as in \fullref{HWknotn}, 
define the map $\Phi\co \Sigma X^{(4)}\lra \Sigma X^{(4)}$ as
\[
\Phi=\tau_{34}-\tau_{12}\tau_{34}+\tau_{14}\tau_{23}+\tau_{13}\tau_{24}
\]
where $\tau$ is the twist map. Then there is a commutative diagram
\[
\xymatrix{ \Omega\Sigma X^{(4)}\ar[d]^{\Omega \Phi}\ar[rr]^{\Omega
W_4} & &\Omega\Sigma X\ar[d]^{H_2}\\
\Omega\Sigma X^{(4)}\ar[rr]^{\Omega W_2(X^{(2)})} && \Omega\Sigma
X^{(2)}}
\]
where $W_n$ denotes the $n$-th fold Whitehead product and $H_2$ the
second James--Hopf map.
\end{thm}

The second application of our method is concerned with exponent
problems in homotopy theory. In general, two types of exponents of a
given space $X$ can be considered. Let $p$ be a prime. The mod~$p$
homotopy exponent of a space $X$ is $p^r$ if that is the least power
of $p$ which annihilates the $p$--torsion component of $\pi_*(X)$.
Denote the mod~$p$ homotopy exponent of $X$ by $\exp(X)=p^r$. A
stronger notion is that of a multiplicative (or $H$-) exponent. If
$Y$ is an $H$--space, then the $p$-th power map is given by the
composite $\ddr{p\co Y}{\Delta}{Y^{\otimes p}}{\mu}{Y}$, where
$\mu$ is the multiplication on~$Y$. The multiplicative exponent of
$Y$ is $p^r$ if that is the least power of $p$ such that $p^r\co
Y\lra Y$ is null homotopic, while $p^{r-1}\co Y\lra Y$ is
essential. We say that $Y$ has no multiplicative exponent if the
$p^r$-th power map on $Y$ is essential for  all $r\in\mathbb{N}$.

A major exponent conjecture was posed by Barratt.
\begin{conj}
Let $f\co \Sigma^2X\lra Y$ be a map of order $p^r$ in
$[\Sigma^2X,Y]$. Then
\[
p^{r+1}\Im \left(f_*\co\pi_*(\Sigma^2 X)\lra\pi_*(Y)\right)=0.
\]
\end{conj}
In particular, if the identity map on $\Sigma^2 X$ is of order $p^r$
in $[\Sigma^2X,\Sigma^2X]$, then the mod~$p$ exponent of $\Sigma^2X$
is $p^{r+1}$, that is,
\[
p^{r+1}\pi_*(\Sigma^2X)=0.
\]
A strong version of the Barratt Conjecture is concerned with the
multiplicative exponent of $\Sigma^2X$ and can be stated as follows.
\begin{conjs}Let $f\co \Sigma^2X\lra Y$ be a map of order $p^r$ in
$[\Sigma^2X,Y]$. Then
\[
\Omega^2f\co \Omega^2\Sigma^2X\lra\Omega^2Y
\]
has order bounded by $p^{r+1}$ in $[\Omega^2\Sigma^2X,\Omega^2Y]$.
\end{conjs}
We start with a map $f\co X\lra \Omega Y$ such that $p^r[f]=0$ in
the group $[X,\Omega Y]$. Let $J(f)\co J(X)\lra \Omega Y$ be its
multiplicative extension to the James construction. Using the
combinatorial description of the group of homotopy classes of maps
from the James filtration $\{J_n(X)\}_{n\geq 0}$ to a loop space, we
estimate the rate of growth of the order of the map $J(f)$.

\begin{thm}
\label{application_Barrattintro} Let $X=\Sigma X'$ be a suspension
and let $f\co X\lra \Omega Y$ be a map such that $p^r[f]=0$ in
the group $[X,\Omega Y]$. Let $J(f)\co J(X)\lra \Omega Y$ be the
canonical multiplicative extension of $f$. Then the following hold.
\begin{enumerate}
 \item The map
\[
J(f)|_{J_n(X)}\co J_n(X)\lra \Omega Y
\]
has order $p^{r+t}$ in $[J_n(X),\Omega Y]$ if $n<p^{t+1}$.
 \item
The composite
\[
\ddr{J_{p^{t+1}}(X)}{J(f)|_{J_{p^{t+1}}(X)}}{\Omega
Y}{p^{r+t}}{\Omega Y}
\]
is homotopic to the composite

\hspace*{-.4cm}$\dddr{J_{p^{t+1}}(X)}{\pinch}{X^{(p^{t+1})}}
{p^{r-1}(\sum_{\tau\in\Sigma_{p^{t+1}-1}}1\wedge
\tau)}{X^{(p^{t+1})}}{\widetilde{W}_{p^{t+1}}}J(X)\stackrel{J(f)}{\lra}{\Omega
Y,}$

\vspace*{.3cm} \noindent where $p^{r+t}\co \Omega Y\lra \Omega Y$
is the $p^{r+t}$-th power map, $\widetilde{W}_n$~is the $n$--fold
Samelson product and

\hspace*{-.4cm}$1\wedge \tau (x_1\wedge \cdots\wedge
x_{p^{t+1}})=x_1\wedge x_{\tau(2)}\wedge\cdots\wedge
x_{\tau(p^{t+1})}\co X^{(p^{t+1})}\lra X^{(p^{t+1})}$

\vspace*{.2cm} \noindent is the map which permutes positions.
\item Let $g=J(f)\circ \widetilde{W}_{p^{t+1}}\circ
(\sum_{\tau\in\Sigma_{p^{t+1}-1}}1\wedge\tau)\circ p^{r-1}\co
X^{(p^{t+1})}\lra\Omega Y$. This is an equivariant map with respect
to the symmetric group action, that is,
\[
g\circ\sigma\simeq g
\]
for any $\sigma\in \Sigma_{p^{t+1}}$.
\end{enumerate}
\end{thm}

It is important to emphasize that parts $(2)$ and $(3)$ of
\fullref{application_Barrattintro} express the first non-trivial
obstruction to the exactness of the Barratt Conjecture in terms of a
computable equivariant map. Some special properties of the trace map
$\Phi=\sum_{\tau\in\Sigma_n}\tau$ allows us to give a more detailed
description of the first obstruction in the case when $X$ is a
two-cell complex (see \fullref{trace} and
\fullref{thm:trace2cell}).

When originally formulated, the Barratt Conjecture did not have any
example which supported it. The first known example that satisfies
the statement of the Barratt Conjecture was the odd primary Moore
space. Neisendorfer~\cite{N}, following Cohen, Moore, and
Neisendorfer's work~\cite{CMN} on the decomposition of the loop
space on the mod~$p^r$ Moore space, showed that $\Omega^2P^n(p^r)$
has multiplicative exponent $p^{r+1}$. When $p$ is an even prime,
Theriault~\cite{Th} showed that when $r\geq 6$ the Moore space
$P^n(2^r)$ has exponent $2^{r+1}$. The existent of bounded exponent
for the mod~$2$ Moore space remains a mystery. Our next goal is to
find a property of $P^{n}(2)$ that will shed some light on the
exponent problem.
\begin{proposition}
\label{mod2Moore} Let $X=P^n(2)$ be the $n$--dimensional mod $2$
Moore space with $n\geq 3$. Then the composite
\[
\ddr{X^{(3)}}{[2]}{X^{(3)}}{\sum_{\sigma\in\Sigma_3}\sigma}{X^{(3)}}
\]
is null homotopic and therefore by
\fullref{application_Barrattintro}~$(2)$, the power map
\[
8|_{J_4(P^n(2))}\co J_4(P^n(2))\lra J(P^n(2))\simeq \Omega
(P^{n+1}(2))
\]
is null homotopic.
\end{proposition}

The paper is organized as follows. In \fullref{sec2}, we translate the
algebraic model of~\cite{GW} into geometry, establishing relations
between combinatorial groups and the groups of natural
transformations of a tensor algebra with the groups of homotopy
classes of maps from the James construction to a loop space.
Applications of this algebraic-geometric model to homotopy theory
are given in Sections~\ref{sec3} and~\ref{sec4}. In \fullref{sec3} we consider the Hopf
invariant of the Whitehead product and prove Theorems~\ref{HWknotn}
and~\ref{HWkn}. The exponent problem, with the emphasis on a rate of
exponent growth for the Barratt Conjecture is treated in \fullref{sec4}.
In particular, in this section we describe some properties of the
$8$-th power map on the mod~$2$ Moore space $P^n(2)$ proving
\fullref{mod2Moore}.

{\bf Acknowledgements}\qua The authors would like to thank Professors
John Berrick, Fred Cohen, Paul Selick and Stephen Theriault for
their helpful suggestions and kind encouragement. The first author
would also like to thank Professor John Berrick and the second
author for making it possible for her to visit the National
University of Singapore for a term and providing her with such a
friendly working atmosphere.
\vspace{-3pt}

\section{Geometrical realisations}
\label{sec2}
\vspace{-3pt}

If $X$ is a connected $CW$ complex, the Bott--Samelson theorem
says that \linebreak $H_*(\Omega\Sigma X)$ is isomorphic as an
algebra to $T(\widetilde{H}_*(X))$, where $T$ denotes the tensor
algebra. We can make $T(\widetilde{H}_*(X))$ into a coalgebra by
requiring that the elements of $\widetilde{H}_*(X)$ are primitive
and then extend to all of $T(\widetilde{H}_*(X))$ via the
multiplication. If $X$ is itself a suspension, then
$T(\widetilde{H}_*(X))$ with this Hopf algebra structure is
isomorphic as a Hopf algebra to $H_*(\Omega\Sigma X)$.
\vspace{-4pt}

\subsubsection*{The James construction}
\vspace{-4pt}
Let $X$ be a topological space with a non-degenerate basepoint $*$
and a compactly generated topology. Then the James construction
$J(X)$ on $X$ is the free topological monoid generated by $X$
subject to the single relation that the basepoint $*$ is the unit.
Combinatorially, the James construction $J(X)$ is obtained from the
disjoint union $\bigsqcup_{k=1}^\infty X^k$ by identifying
$(x_1,\ldots ,x_i,\ldots ,x_k)$ with $(x_1,\ldots ,\hat{x}_i,\ldots
,x_k)$ if $x_i=*$. Non-unit points of $J(X)$ can thus be thought of
as words $(x_1,\ldots ,x_k)$ of length $k$, with no $x_i$ being the
unit. Let $q_n\co X^n\lra J_n(X)$ be the quotient map. The James
filtration $\{ J_n(X)\}_{n\geq 0}$ with $J_0(X)=*$ and $J_1(X)=X$ is
induced by the word-length filtration. It follows readily that
$J_n(X)/J_{n-1}(X)$ is homeomorphic to $X^{(n)}$. The fundamental
properties of $J(X)$ are as follows:
\begin{enumerate}
 \item if $X$ is path connected, then $J(X)$ is (weak) homotopy equivalent
 to $\Omega\Sigma X$;
 \item the quotient $\Sigma q_n\co \Sigma X^n\lra \Sigma J_n(X)$ has
 a functorial cross-section;
 \item the inclusion $\Sigma J_{n-1}(X)\lra \Sigma J_n(X)$ has a
 functorial retraction;
 \item there is a functorial decomposition
 \begin{equation}
 \label{stablesplitting}
 \Sigma J_n(X)\simeq \bigvee^n_{j=1}\Sigma X^{(j)}
 \end{equation}
 for $0\leq n\leq\infty$.
\end{enumerate}
\vspace{-4pt}
Recall that a space $X$ has a weak $LS$--category less than $k$ if
the reduced diagonal $\bar\Delta_k\co X\lra X^{(k)}$ is null
homotopic. The spaces that we will consider in the following
sections will be path connected and will have a weak $LS$--category
less than~2, that is, the reduced diagonal $\bar\Delta\co X\lra
X\wedge X$ is null homotopic (for example, if $X$ is a co-$H$
space).
\vspace{-4pt}

\subsubsection*{Geometrical Cohen groups}
\vspace{-4pt}
Let $X$ be a pointed path connected space, $X^n$ the $n$--fold self
Cartesian product of $X$ and $X^{(k)}$ the $k$--fold self smash
product of $X$. Let $J(X)$ be the James construction of $X$. The
following groups were introduced by Cohen in~\cite{Cohen} and
Wu~\cite{Wu}.
\vspace{-4pt}
\begin{defn}
Let $K_n(X)$ denote the subgroup of $[X^n, J(X)]$ generated by the
homotopy classes $x_i$ for $1\leq i\leq n$, where $x_i$ is
represented by the composite
\[
\ddr{X^n}{p_i}{X}{E}{J(X)}
\]
where $p_i\co X^n\lra X$ is the $i$-th coordinate projection
given by
\[
p_i(x_1,x_2,\cdots, x_n)=x_i,
\]
and $E\co X\lra J(X)$ is the canonical inclusion.
\end{defn}
\vspace{-4pt}

\begin{proposition}[Cohen~\cite{Cohen}]\label{prop:geo_relations}
Let $X$ be a path connected space with weak $LS$--category less
than~2. Then, in the group $[X^n,JX]$, the following identities
hold:
\begin{enumerate}
\item\strutt $[[x_{i_1},x_{i_2},\cdots,x_{i_k}]=1$\quad
if $i_s=i_t$ for some $1\leq s<t\leq k$,\\\strutt where\ \
$[[a_1,a_2,\cdots,a_l]\negthinspace=\negthinspace[\cdots[a_1,a_2],\cdots\negthinspace,],a_l]$\ \ with $[x,y]=x^{-1}y^{-1}xy$; 

\item\strutt $[[x_{i_1}^{n_1},x_{i_2}^{n_2},\cdots,x_{i_k}^{n_l}]=
[[x_{i_1},x_{i_2},\cdots,x_{i_k}]^{n_1n_2\cdots n_l}.$
\end{enumerate}
\end{proposition}
Recall that the Cohen group $K_n$ is defined combinatorially as
follows. The Cohen group $K_n(x_1,x_2,\ldots,x_n)$ is the quotient
group of the free group $F_n$ of rank $n$ generated by
$x_1,x_2,\ldots,x_n$ modulo the relations
\begin{enumerate}
 \item $[[x_{i_1},x_{i_2},x_{i_3},\cdots,x_{i_l}]=1$\quad if $i_s=i_t$
 for some $1\leq s,t\leq l$;
 \item
 $[[x_{i_1}^{n_1},x_{i_2}^{n_2},x_{i_3}^{n_3},\cdots,x_{i_l}^{n_l}]=
 [[x_{i_1},x_{i_2},x_{i_3},\cdots,x_{i_l}]^{n_1n_2\cdots n_l}$.
\end{enumerate}
\begin{proposition}
\label{prop:georealK_n} Let $X$ be a path connected space with weak
$LS$--category less than~$2$. Then there is a homomorphism
\[
e_X\co K_n \lra K_n(X)\subseteq[X^n, J(X)]
\]
given on any generator $x_i$ of $K_n$ by $e_X(x_i)=x_{i}.$
\end{proposition}
\begin{proof}
The existence of the homomorphism $e_X$ follows readily from the
definition of $K_n$ and \fullref{prop:geo_relations}.
\end{proof}

This homomorphism can be generalized in the following way.
\begin{corollary}
Let $X$ be a path connected space with weak $LS$--category less
than~$1$, $M$ a path connected topological monoid and $f\co X\lra
M$ a pointed map. Then there is homomorphism
\begin{gather*}
e_f\co K_n \lra\,[X^n, M]
\\
\tag*{\hbox{given by}}
e_f(x_i)=J(f)_*(x_{i}),
\end{gather*}
where $J(f)\co J(X^{(k)})\lra M$ is a homomorphism of topological
monoids such that $J(f)|_{X}=f$ and $J(f)_*\co[X^n,
J(X)]\lra\,[X^n,M]$ is induced by the map $J(f)$.
\end{corollary}

\begin{proposition}[Cohen~\cite{Cohen}]
Let $X$ be a path connected space with weak $LS$--category less
than~$1$, $M$ a path-connected topological monoid and $f\co X\lra
M$ a pointed map. Suppose that the $q$-th power $[f]^q=1$ in the
group $[X,M]$. Then the homomorphism $\theta_f\co
K_n\lra\,[X^n,M]$ factors through the quotient group $K_n^{\Z/q}$.
\end{proposition}

\subsubsection*{Generalization to $K_n(k)(X)$}

\begin{defn}
\hspace*{-1mm}Let $k\negthinspace\leq\negthinspace n$. Define
$K_n(k)(X\negthinspace)$ to be the subgroup of $[X^n\negthinspace,
J(X^{(k)}\negthinspace)]$ generated by the homotopy classes
$\{x_{i_1}|x_{i_2}|\cdots|x_{i_k}\}$ for $1\leq i_j\leq n$ and
$1\leq j\leq k$, where $\{x_{i_1}|x_{i_2}|\cdots|x_{i_k}\}$ is
represented by the composite
\[
\ddr{X^n}{p_{i_1\cdots i_k}}{X^{(k)}}{E}{J(X^{(k)})},
\]
where $E\co X^{(k)}\lra J(X^{(k)})$ is the canonical inclusion
and $p_{i_1\cdots i_k}\co X^n\lra X^{(k)}$ is given by
\[
p_{i_1\cdots i_k}(x_1,x_2,\cdots, x_n)=x_{i_1}\wedge\cdots\wedge
x_{i_k}.
\]
\end{defn}
Notice that in the case $k=1$, the group $K_n(1)(X)$ is the Cohen
group $K_n(X)$.

\begin{proposition}[Lemma 2.2, ~\cite{Wu}]\label{prop:geo_relationK(k)}
Let $X$ be a path connected space with weak $LS$--category less
than~2. Then, in the group $[X^n,J(X^{(k)})]$, the following
identities hold:
\begin{enumerate}
\item\strutt $\{x_{i_1}|x_{i_2}|\cdots|x_{i_k}\}=1$\quad
 if $i_s=i_t$ for some $1\leq s<t\leq k$;
\item
\strutt$[[\{x_{i_1}|x_{i_2}|\cdots|x_{i_k}\},\{x_{i_{k+1}}|x_{i_{k+2}}
 |\cdots|x_{i_{2k}}\},\negthinspace\cdots$\\
\strutt\hspace{2in}$\cdots\negthinspace,\{x_{i_{(l-1)k+1}}|x_{i_{(l-1)k+2}}|
\cdots|x_{i_{lk}}\}]=1$\\
\strutt if $i_s=i_t$ for some $1\leq s<t\leq kl$, where
 $[[a_1,a_2,\cdots\negthinspace,a_l]=[\cdots[a_1,a_2],\cdots\negthinspace,],a_l]$ with
 $[x,y]=x^{-1}y^{-1}xy$;

\item\strutt$[[\{x_{i_1}|x_{i_2}|\negthinspace\cdots\negthinspace|x_{i_k}\}^{n_1},\{x_{i_{k+1}}|x_{i_{k+2}}
 |\negthinspace\cdots\negthinspace|x_{i_{2k}}\}^{n_2},\negthinspace\cdots$\\
\strutt\hspace{2in}$\cdots\negthinspace,\{x_{i_{(l-1)k+1}}|x_{i_{(l-1)k+2}}|\negthinspace\cdots\negthinspace|x_{i_{lk}}\}^{n_l}]=$\\
$[[\{x_{i_1}|x_{i_2}|\negthinspace\cdots\negthinspace|x_{i_k}\},\{x_{i_{k+1}}|x_{i_{k+2}}|\negthinspace\cdots\negthinspace|
 x_{i_{2k}}\},\negthinspace\cdots$\\
\strutt\hspace{2in}$\cdots\negthinspace,\{x_{i_{(l-1)k+1}}|x_{i_{(l-1)k+2}}|\negthinspace\cdots\negthinspace|x_{i_{lk}}\}]^{n_1n_2\cdots
 n_l}_.$
\end{enumerate}
\end{proposition}

\begin{rem}
Relation $(3)$ follows from relation $(2)$ (see for
example~\cite{Cohen}).
\end{rem}

Recall from~\cite{GW} that the group $K_n(k)$ is defined
combinatorially as follows:
\begin{enumerate}
\item
 generators are the words $\{x_{i_1}|x_{i_2}|\cdots|x_{i_k}\}$ with $1\leq i_j\leq n$ for $1\leq j\leq
 k$;

\item
 relations are given by identities $(1)-(3)$ in
 \fullref{prop:geo_relationK(k)}.\newline

\noindent Let $q$ be an integer. The group $K_n^{\Z/q}(k)$ is the
quotient group of $K_n(k)$ modulo the following additional
relations:
\item
 $\{x_{i_1}|x_{i_2}|\cdots|x_{i_k}\}^q=1$ for each generator $\{x_{i_1}|x_{i_2}|\cdots|x_{i_k}\}$.
\end{enumerate}
\begin{proposition}
\label{prop:georealK_n(k)} Let $X$ be a path connected space with
weak $LS$--category less than~$2$. Then there is a homomorphism
\[
e_X\co K_n(k) \lra[X^n, J(X^{(k)})]
\]
given on any generator $\{x_{i_1}|x_{i_2}|\cdots|x_{i_k}\}$ of
$K_n(k)$ by
\[
e_X(\{x_{i_1}|x_{i_2}|\cdots|x_{i_k}\})=\{x_{i_1}|x_{i_2}|\cdots|x_{i_k}\}.
\]
\end{proposition}
\begin{proof}
The existence of the homomorphism $e_X$ follows immediately from the
definition of $K_n(k)$ and \fullref{prop:geo_relationK(k)}.
\end{proof}
\begin{corollary}
Let $X$ be a path connected space with weak $LS$--category less
than~$2$, $M$ be a path connected topological monoid and $f\co
X^{(k)}\lra M$ a pointed map. By the definition of $K_n(k)$, the
homomorphism
\begin{gather*}
e_f\co K_n(k)\lra[X^n, M],
\\
\tag*{\hbox{given by}}
e_f(\{x_{i_1}|x_{i_2}|\cdots|x_{i_k}\})=J(f)_*(\{x_{i_1}|x_{i_2}|\cdots|x_{i_k}\}),
\end{gather*}
is well-defined, where $J(f)\co J(X^{(k)})\lra M$ is a
homomorphism of topological monoids such that $J(f)|_{X^{(k)}}=f$
and $J(f)_*\co[X^n, J(X^{(k)})]\lra[X^n,M]$ is induced by the map
$J(f)$.
\end{corollary}

\begin{lem}
Let $X$  be a path connected space with weak $LS$--category less
than~2, $M$ a path connected topological monoid and $f\co
X^{(k)}\lra M$ a pointed map. Suppose that the $q$-th power
$[f]^q=1$ in the group $[X^{(k)},M]$. Then the homomorphism
$e_f\co K_n(k)\lra[X^n,M]$ factors through the quotient group
$K_n^{\Z/q}(k)$.
\end{lem}
\vspace{-4pt}
\begin{proof}
The image under $e_f$ of the element
$\{x_{i_1}|x_{i_2}|\cdots|x_{i_k}\}^q$ is represented by the
composite
\[
\ddddr{X^n}{p_{i_1\cdots i_k}}{X^{(k)}}{E}{J(X^{(k)})}{J(f)}{
J(M)}{q}{J(M)}\lra M.
\]
As $[f]^q=1$, the composite is null homotopic. Thus the assertion
follows.
\end{proof}
\vspace{-4pt}

\subsubsection*{Generalization to $\H_n(X)$}
\vspace{-4pt}

By the suspension splitting theorem \ref{stablesplitting}, the
inclusion
\[
J_{n-1}(X)\lra J_n(X)
\]
induces a tower of group epimorphisms
\[
[J(X),\Omega Y] \lra\cdots\lra\,[J_n(X),\Omega Y] \lra\cdots \lra\,
[X,\Omega Y],
\]
and there is a group isomorphism between $[J(X),\Omega Y]$ and the
inverse limit
\[
[J(X),\Omega Y]\cong\lim_n [J_n(X),\Omega Y].
\]
Let $q_n\co X^n\lra J_n(X)$ be the quotient map. Then by the
suspension splitting theorem \ref{stablesplitting}, there is a group
monomorphism
\[
q^*_n\co [J_n(X), J(X)]\lra\,[X^n, J(X)]
\]
for each $n$.
\begin{defn}
Let $\H_n(X)$ be the subgroup of $[J_n(X), J(X)]$ defined by
\[
[J_n(X),J(X)]\cap K_n(X).
\]
\end{defn}
There is an equivalent definition for the groups $\H_n(X)$ which is
more suitable for generalization. Recall that $J_n(X)$ is the
coequalizer of the inclusions $i_j\co X^{n-1}\lra X^n$ for $1\leq
j\leq n$. These inclusions induce projections
\[
i^*_j\co [X^n, Y]\lra\,[X^{n-1},Y]
\]
for $1\leq j\leq n$ and any space $Y$. Thus there are projections
$d_j\co K_n(X)\lra K_{n-1}(X)$ given by
\[
d_j(x_i)=\left\{\begin{array}{lll}
x_i & \text{for} & i<j\\
1 & \text{for} & i=j\\
x_{i-1} & \text{for} & i>j\qquad\qquad\text{for $1\leq j\leq n$.}
 \end{array}\right.
\]
\begin{lem}
The diagram of group homomorphisms
\begin{equation}\label{dgm:K_nK_n(X)}
\xymatrix{
K_n \ar[r]^-{e_X} \ar@{->>}[d]^{dj}& K_n(X)\ar@{->>}[d]^{dj}\\
K_{n-1}\ar[r]^-{e_X} & K_{n-1}(X)}
\end{equation}
commutes for every $1\leq j\leq n$.
\end{lem}
\vspace{6pt}

\begin{proof}
The proof follows immediately from the definitions.
\end{proof}
\vspace{6pt}

Define the group $\H_n(X)$ to be the equalizer of homomorphisms
$d_j$ for $1\leq j\leq n$. From the above discussion, it follows
readily that the two definitions of $\H_n(X)$ are equivalent. Now
from the second definition of $\H_n(X)$, as $d_i\mid_{\H_n(X)}=
d_j\mid_{\H_n(X)}$ for $1\leq i,j\leq n$, there are homomorphisms
$p_n\co\H_n(X)\lra\H_{n-1}(X)$ such that the diagram
\[
\xymatrix{
\H_n(X)\ar[r]\ar[d]^{p_n} & K_n(X)\ar@{->>}[d]^{d_j}\\
\H_{n-1}(X)\ar[r]& K_{n-1}(X)}
\]
commutes for $1\leq j\leq n$. Selick and Wu~\cite{SW} proved that
there is a progroup
\[
\H_{\infty}(X) \lra\cdots\lra \H_{n}(X)\lra\cdots \lra\H_{0}(X)
\]
where $\H_{\infty}(X)$ is given by the inverse limit
\[
\H_{\infty}(X)=\lim_n \H_n(X).
\]
\medskip

In~\cite{GW} we recalled the combinatorial group $\H_n$ as the
equalizer of the projections $p_j\co K_n\lra K_{n-1}$ for $1\leq
j\leq n$. The definition of $\H_n$ in this setting goes back to
Cohen~\cite{Cohen}.
\medskip

\begin{proposition}
\label{prop:georealH_n} Let $X$ be a path connected space with weak
$LS$--category less than~$2$.

Then there is a homomorphism
\[
e_X\co \H_n \lra\H_n(X)\subseteq [J_n(X), J(X)]
\]
such that the diagram
\[
\xymatrix{\H_n \ar[r]^-{e_X}\ar[d] & [J_n(X), J(X)]\ar[d]\\
K_n \ar[r]^-{e_X} & [X^n, J(X)]}
\]
commutes.
\end{proposition}
\begin{proof}
Notice that the group $[J_n(X), J(X)]$ can be thought of as the
equalizer of the projections $i_j^*\co [X^n, J(X)]\lra\,[X^{n-1},
J(X)]$. The existence of the homomorphism $e_X\co \H_n
\lra\,[J_n(X), J(X)]$ follows readily from the existence of
$e_X\co K_n \lra K_n(X)\subseteq [X^n, J(X)]$, the fact that the
groups $\H_n$ and $\H_n(X)$ are the equalizers of appropriate
projections on $K_n$ and $K_n(X)$, respectively , and that
diagram~\eqref{dgm:K_nK_n(X)} commutes.
\end{proof}
\begin{corollary}
Let $X$ be a path connected space with weak $LS$--category less
than~$2$, $M$ be a path connected topological monoid and $f\co
X\lra M$ a pointed map. Then there is a homomorphism
\[
e_f\co\H_n \lra\,[J_n(X), M]
\]
given by $e_f= J(f)_*\circ e_X$.
\end{corollary}

\begin{corollary}\label{cor:H_n^qtoM}
Let $X$  be a path connected space with weak $LS$--category less
than~$2$, $M$ a path connected topological monoid and $f\co X\lra
M$ a pointed map. Suppose that the $q$-th power $[f]^q=1$ in the
group $[X,M]$. Then the homomorphism $e_f\co
\H_n\lra\,[J_n(X),M]$ factors through the quotient group
$\H_n^{\Z/q}$.
\end{corollary}

\subsubsection*{Generalization to $\H_n^{(l)}(X)$}
Further on, we want to find a geometrical analogue of the
combinatorial group $\H^{(l)}_n$ defined in~\cite[Definition
2.7]{GW}.
\begin{proposition}
The group $[J_n(X^{(l)}),J(X)]$ is the equalizer of the restriction
of the projections $D_j\co K_{ln}(X)\lra K_{l(n-1)}(X)$, that is,
\[
\xymatrix{
K_{ln}(X)\bigcap\left(\bigcap_{s=0}^{n-1}\left(\bigcap^l_{j=1}\Ker
i^*_{sl+1}\right) \right)\ar@<-5pt>[d]_{D_0}^{..}\ar@<5pt>[d]_{.}^{D_{n-1}}\\
K_{l(n-1)}(X)\bigcap\left(\bigcap_{s=0}^{n-2}
\left(\bigcap^l_{j=1}\Ker i^*_{sl+1}\right)\right)}
\]
where the projection $D_j$ is given by
\[
D_j(x_i)=\left\{\begin{array}{lll}
x_i & \text{for} & i<lj\\
1 & \text{for} & lj+1\leq i\leq (l+1)j\\
x_{i-1} & \text{for} & (l+1)j<i
\end{array}\right.
\] for $0\leq j\leq n-1$.
\end{proposition}
\begin{proof}
Notice that the composite
\[
X^{ln}\lra\underbrace{X^{(l)}\times\cdots\times X^{(l)}}_n \lra
J_n(X^{(l)})
\]
induces a group monomorphism
\[
[J_n(X^{(l)}),J(X)]\lra\,[\underbrace{X^{(l)}\times\cdots\times
X^{(l)}}_n, J(X)]\lra[X^{ln}, J(X)].
\]
Let $d_j$ denote the following projection
\[
\underbrace{X^{(l)}\times\cdots\times X^{(l)}}_n\lra
\underbrace{X^{(l)}\times\cdots\times \widehat{X^{(l)}}\cdots\times
X^{(l)}}_{n-1}
\]
for $1\leq j\leq n$.

Let $i_j\co X^l\lra X^{l+1}$ denote the usual coordinate
inclusion
\[
i_j(x_1,\ldots,x_l)=(x_1,\ldots,x_{j-1},*, x_j,\ldots x_l).
\]
Notice that if
\[
\alpha\in [X^{ln},J(X)]\cap [\underbrace{X^{(l)}\times\cdots\times
X^{(l)}}_n, J(X)],
\]
then
\[
\alpha\in \bigcap_{s=0}^{n-1}\left(\bigcap^l_{j=1}\Ker
i^*_{sl+j}\right).
\]
Now it is clear that $[J_n(X^{(l)}),J(X)]$ is
\[
\mathrm{eq}(D_0,\ldots,D_{n-1})\bigcap\left(\bigcap_{s=0}^{n-1}\left(\bigcap^l_{j=1}\Ker
i^*_{sl+j}\right)\right)
\]
where $\mathrm{eq}(f_1,\ldots,f_n)$ stands for the equalizer of maps
$f_1,\ldots,f_n$.
\end{proof}
\begin{defn}\label{def:Hln(X)}
Define the group
$\H^{(l)}_n\negthinspace(\negthinspace X\negthinspace)$ to be the
subgroup of $[J_n(X^{(l)}\negthinspace), J(X\negthinspace)]$ given
by the equalizer of the projections
\[
\xymatrix{ K_{ln}(X)\bigcap\left(\bigcap_{s=0}^{n-1}\left(
\bigcap^l_{j=1}\Ker
i^*_{sl+1}\right)\right)\ar@<-5pt>[d]_{D_0}^{..}\ar@<5pt>[d]_{.}^{D_{n-1}}\\
K_{l(n-1)}(X)\bigcap\left(\bigcap_{s=0}^{n-2}\left(\bigcap^l_{j=1}\Ker
i^*_{sl+1} \right)\right).}
\]
\end{defn}
From the definition of $\H_n^{(l)}(X)$, as
$D_i\mid_{\H_n^{(l)}(X)}=D_j\mid_{\H_n^{(l)}(X)}$ for every $1\leq
i,j\leq n$, there is a homomorphism
$p_n\co\H_n^{(l)}(X)\lra\H_{n-1}^{(l)}(X)$ such that the diagram
\begin{equation}\label{dgmK_nH_n}
\xymatrix{ \H_n^{(l)}(X)\ar[r]\ar[d]^{p_n}&
K_{ln}(X)\ar@{->>}[d]^{D_j}\\
\H_{n-1}^{(l)}(X)\ar[r]& K_{l(n-1)}(X)}
\end{equation}
commutes for $0\leq j\leq n-1$.
\begin{lem}
There is a progroup
\[
\H^{(l)}(X) \lra\cdots\lra \H_{n}^{(l)}(X)\lra\cdots
\lra\H_{0}^{(l)}(X)
\]
where $\H^{(l)}(X)$ is given by the inverse limit
\[
\H^{(l)}(X)=\lim_{p_n} \H_n^{(l)}(X).
\]
\end{lem}
\begin{proof}
The map $p_n\co\H_n^{(l)}(X)\lra\H_{n-1}^{(l)}(X)$ is obtained as
the map of equalizers induced by the epimorphisms $D_j$. Thus $p_n$
is an epimorphism for every $n>0$.
\end{proof}
\begin{proposition}
\label{prop:georealH_n^l} For every $n>0$, there is a homomorphism
\[
e_X\co \H^{(l)}_n\lra\H^{(l)}_n(X)\subseteq [J_n(X^{(l)}),J(X)]
\]
such that the diagram
\[
\xymatrix{ \H^{(l)}_n \ar[r]^-{e_X}\ar[d]  &
[J_n(X^{(l)}),J(X)]\ar[d]\\
K_{ln}\ar[r]^-{e_X} & [X^{ln},J(X)]}
\]
commutes.
\end{proposition}
\begin{proof}
The statement follows immediately from the definition of the
groups $\H_n^{(l)}$ and $\negthinspace\H_n^{(l)}(X\negthinspace)$
and the existence of the homomorphism $e_X\co
K_{ln}\negthinspace\negthinspace\to\negthinspace[X^{ln}\negthinspace,J(\negthinspace
X\negthinspace)]$.
\end{proof}

\begin{corollary}
Let $X$ be a path connected space with weak $LS$--category less
than~$1$, $M$ be a path connected topological monoid and $f\co
X\lra M$ a pointed map. Then there is a homomorphism
\[
e_f\co\H_n^{(l)} \lra[J_n(X^{(l)}), M]
\]
given by $e_f= J(f)_*\circ e_X$.
\end{corollary}

\begin{corollary}
Let $X$  be a path connected space with weak $LS$--category less
than~2, $M$ a path connected topological monoid and $f\co X\lra
M$ a pointed map. Suppose that the $q$-th power $[f]^q=1$ in the
group $[X,M]$. Then the homomorphism $e_f\co
\H^{(l)}_n\lra[J_n(X^{(l)}),M]$ factors through the quotient group
$^{\Z/q}\H^{(l)}_n$.
\end{corollary}

\subsubsection*{Generalization to $\H_n^{(l)(k)}(X)$}

\begin{defn}
Define the group $\H^{(l)(k)}_n(X)$ to be the subgroup of the
group\break $[J_n(X^{(l)}), J(X^{(k)})]$ given by the equalizer of the
projections
\[
\xymatrix{
K_{ln}(k)(X)\bigcap\left(\bigcap_{s=0}^{n-1}\left(\bigcap^l_{j=1}\Ker
i^*_{sl+1}\right)\right)
\ar@<-5pt>[d]_{D_0}^{..}\ar@<5pt>[d]_{.}^{D_{n-1}}\\
K_{l(n-1)}(k)(X)\bigcap\left(\bigcap_{s=0}^{n-2}\left(\bigcap^l_{j=1}\Ker
i^*_{sl+1}\right)\right).}
\]
where the homomorphisms $D_i$ are induced by those from
\fullref{def:Hln(X)}.
\end{defn}
From the definition of $\H_n^{(l)(k)}(X)$, as
$D_i\mid_{\H_n^{(l)(k)}(X)}=D_j\mid_{\H_n^{(l)(k)}(X)}$ for every
$1\leq i,j\leq n$, there is a homomorphism $p_n\co
\H^{(l)(k)}_n(X)\lra \H^{(l)(k)}_{n-1}(X)$ such that the following
diagram
\[
\xymatrix{ \H_n^{(l)(k)}(X)\ar[r] \ar[d]^{p_n}&
K_{ln}(k)(X)\ar@{->>}[d]^{D_j}\\
\H_{n-1}^{(l)(k)}(X)\ar[r]& K_{l(n-1)}(k)(X)}
\]
commutes for $0\leq j\leq n-1$.
\begin{lem}
There is a progroup
\[
\H^{(l)(k)}(X) \lra\cdots\lra \H_{n}^{(l)(k)}(X)\lra\cdots
\lra\H_{0}^{(l)(k)}(X)
\]
where $\H^{(l)(k)}(X)$ is given by the inverse limit
\[
\H^{(l)(k)}(X)=\lim_{p_n} \H_n^{(l)(k)}(X).
\]
\end{lem}
\begin{proof}
The map $p_n\co\H_n^{(l)(k)}(X)\lra\H_{n-1}^{(l)(k)}(X)$ is
obtained as the map of equalizers induced by the epimorphisms $D_j$.
Thus $p_n$ is an epimorphism for every $n>0$.
\end{proof}
The combinatorial group $\H_n^{(l)(k)}$ is defined in the
predecessor paper~\cite[Definition~2.9]{GW}.
\begin{proposition}
\label{prop:georealH_n^lk} There is a homomorphism
\[
e_X\co \H^{(l)(k)}_n\lra\H^{(l)(k)}_n(X)\subseteq
[J_n(X^{(l)}),J(X^{(k)})]
\]
such that the diagram
\[
\xymatrix{ \H^{(l)(k)}_n \ar[r]^-{e_X}\ar[d]  &
[J_n(X^{(l)}),J(X^{(k)})]\ar[d]\\
K_{ln}(k)\ar[r]^-{e_X} & [X^{ln},J(X^{(k)})]}
\]
commutes.
\end{proposition}
\begin{proof}
The statement follows immediately from the definition of the group
$\H^{(l)(k)}\negthinspace(X\negthinspace)$ and the existence of
the homomorphism $e_X\co
K_{ln}(k)\negthinspace\to\negthinspace[X^{ln}\negthinspace,J(X^{(k)}\negthinspace)]$.
\end{proof}
\begin{corollary}
Let $X$ be a path connected space with weak $LS$--category less
than~2, $M$ be a path connected topological monoid and $f\co
X^{(k)}\lra M$ a pointed map. Then there is a homomorphism
\[
e_f\co\H_n^{(l)(k)} \lra\,[J_n(X^{(l)}), M]
\]
given by $e_f= J(f)_*\circ e_X$.
\end{corollary}

\begin{corollary}
Let $X$ be a path connected space with weak $LS$--category less
than~2, $M$ a path connected topological monoid and $f\co
X^{(k)}\lra M$ a pointed map. Suppose that the $q$-th power
$[f]^q=1$ in the group $[X^{(k)},M]$. Then the homomorphism
$e_f\co \H^{(l)(k)}_n\lra\,[J_n(X^{(l)}),M]$ factors through the
quotient group $^{\Z/q}\H^{(l)(k)}_n$.
\end{corollary}

\section{Application: Whitehead products and James--Hopf maps}
\label{sec3}

\subsubsection*{James--Hopf maps}

Let $X$ be a pointed space. The James--Hopf map
\[
H_k\co J(X)\lra J(X^{(k)})
\]
is combinatorially defined by
\[
H_k(x_1x_2\ldots x_n)=\prod_{1\leq i_1<i_2<\cdots <i_k\leq
n}(x_{i_1}x_{i_2}\ldots x_{i_k})
\]
with right lexicographical order in the product.

The James--Hopf map can be also defined more geometrically. Let $X$
have a homotopy type of a $CW$ complex. Then using the simple
combinatorial structure of the James construction $J(X)$, there can
be made a preferred choice (inductively on the James filtration
$J_n(X)$) of the homotopy equivalence
\[
\theta\co\Sigma(J(X))\simeq\Sigma\Omega\Sigma X \lra
\bigvee_{i=0}^\infty \Sigma X^{(i)}.
\]
Furthermore, consider the adjoint of $\theta$
\[
\overline{\theta}\co\Omega\Sigma X\lra \Omega\bigl(
\bigvee_{i=0}^\infty \Sigma X^{(i)}\bigr),
\]
and the pinch map
\[
q_k\co \bigvee_{i=0}^{\infty} X^{(i)}\lra X^{(k)}
\]
which sends $X^{(i)}$ to the base point if $i\neq k$ and which is
the identity when restricted to $X^{(k)}$. Now, the $k$-th
James--Hopf map
\[
H_k\co \Omega\Sigma X\lra \Omega\Sigma X^{(k)}
\]
is defined as the composite
$\Omega\Sigma(q_k)\circ\overline{\theta}$.

In this paper we shall take advantage of the combinatorial
definition of the James--Hopf map.

For each $n$, let ${H_k}_*\co [X^n, J(X)]\lra\,[X^n,J(X^{(k)})]$
be the function induced by the map $H_k$. Notice that ${H_k}_*$ is
not a homomorphism of groups if $n\geq k>1$ and $X$ is a
non-contractible suspension. We want to study the map $H_{k*}$ using
the combinatorial methods described in \fullref{sec1}, namely, the
relation between groups $K_n$ and $K_n(k)$ and the natural maps from
$X^n$ to $J(X)$, $J(X^{(k)})$, respectively. Therefore we proceed by
defining the combinatorial analogue
\[
H_k\co K_n=K_n(1)\lra K_n(k)
\]
of $H_{k*}$.

\begin{defn}
The function $H_k\co K_n=K_n(1)\lra K_n(k)$ is defined by setting
\[
H_k(x^{n_1}_{i_1}x^{n_2}_{i_2}\cdots x^{n_l}_{i_l})=\prod_{1\leq
j_1<\cdots<j_k\leq
l}\{x_{i_{j_1}}|x_{i_{j_2}}|\cdots|x_{i_{j_k}}\}^{n_{j_1}n_{j_2}\cdots
n_{j_k}}
\]
with right lexicographical order, for any word
$x^{n_1}_{i_1}x^{n_2}_{i_2}\cdots x^{n_l}_{i_l}\in K_n$.
\end{defn}

\begin{proposition}{\rm\cite[Lemma 2.3]{Wu}}\label{proposition5.2}\qua
The function $H_k\co K_n\lra K_n(k)$ is well-defined.
Furthermore, there is a commutative diagram
\[
\xymatrix{ K_n\ar[r]^-{\theta}\ar[d]^{H_k}& K_n(X)\ar[r]&
[X^n,J(X)]\ar[d]^{{H_k}_*}\\
K_n(k)\ar[r]^-{\theta} &K_n(k)(X)\ar[r]& [X^n,J(X^{(k)})]}
\]
for any suspension $X$.
\end{proposition}

The algebraic analogue of the James--Hopf map $H_k\co T(V)\lra
T(V^{\otimes k})$ is the functorial coalgebra map induced from
\[
H_{k*}\co H_*(\Omega\Sigma X)\lra H_*(\Omega\Sigma (X^{(k)})).
\]
\begin{lem}
The function $H_k\co K_{nt}\lra K_{nt}(k)$ induces a function
\[
h_k\co \H^{(n)}_t\lra\H^{(n),(k)}_t
\]
for any $1\leq t\leq\infty$. Moreover there is a commutative diagram
\begin{equation}
\label{Hdgrm} \xymatrix{ [J(X^{(n)}), J(X)] \ar[r]^-{H_{k*}} &
[J(X^{(n)}), J(X^{(k)})]\\
\H^{(n)}\ar[u]\ar[d]^{\cong}\ar[r]^-{h_k} & \H^{(n),(k)}\ar[u]\ar[d]^{\cong}\\
\Coalg(T(V^{\otimes n}), T(V))\ar[r]^-{H_{k*}} & \Coalg(T(V^{\otimes
n}), T(V^{\otimes k})).}
\end{equation}
\end{lem}
\begin{proof}
The first statement follows from the direct computation and the fact
that $\H^{(n)}_t$ and $\H^{(n),(k)}_t$ are given by certain
equalizers~\cite[Definitions 2.7, 2.9]{GW}.

The top square of the diagram commutes by
\fullref{proposition5.2}. In paper~\cite{GW}, we established
the progroup isomorphisms $e\co\H^{(n)}\lra \Coalg(T(V^{\otimes
n}), T(V))$ and $e\co\H^{(n),(k)}\lra \Coalg(T(V^{\otimes n}),
T(V^{\otimes k}))$. By applying the homology functor, it follows
that
\[
\xymatrix{ \H^{(n)}\ar[d]\ar[r]^-{h_k} & \H^{(n),(k)}\ar[d]\\
[H_*(J(X^{(n)})), H_*(J(X))] \ar[r]^-{H_{k*}} & [H_*(J(X^{(n)})),
H_*(J(X^{(k)}))]}
\]
and hence the bottom square commutes by letting $X$ run through
wedges of the $2$--sphere.
\end{proof}

The main objective of this section is the study of the Hopf
invariant of a Whitehead product. Let $W_n\co\Sigma X^{(n)}\lra
\Sigma X$ denote the $n$ fold Whitehead product on $X$. Our main
result is as follows.
\begin{thm}\label{Thm3.4}
Let $X$ be a pointed space  with the null homotopic reduced diagonal
$\bar{\Delta}\co X\lra X\wedge X$. Then for $n>k$,
\[
\ddr{J(X^{(n)})}{\Omega W_n}{J(X)}{H_k}{J(X^{(k)})}
\]
is null homotopic if $k$ does not divide $n$.
\end{thm}
\begin{proof}
We start by calculating the algebraic analogue of the $k$-th Hopf
invariant of the $n$--fold Whitehead product and then transfer the
result to topology. The geometrical map $H_k\circ\Omega W_n$ induces
the map
\[
H_k\circ T(\beta_n)\co T(V^{\otimes n})\lra T(V^{\otimes k}).
\]
Wu~\cite{Wu} proved that the map $H_k\circ\Omega\omega_n$ is a loop
map. Therefore the induced map $H_k\circ T(\beta_n)$ is an algebra
map which is determined by its values on $V^{\otimes n}$. If $n$ in
not a multiple of $k$, then because of dimensional reasons
$H_k\circ\beta_n=0$, that is, $H_k\circ T(\beta_n)\mid_{V^{\otimes
n}}=0$. That shows that $H_k\circ T(\beta_n)=0$. Now applying
diagram~\eqref{Hdgrm} proves the theorem.
\end{proof}

\begin{thm}
For $X$ as in \fullref{Thm3.4}, define the map $\Phi\co 
\Sigma X^{(4)}\lra \Sigma X^{(4)}$ as
\[
\Phi=\tau_{34}-\tau_{12}\tau_{34}+\tau_{14}\tau_{23}+\tau_{13}\tau_{24}.
\]
Then there is a commutative diagram
\[
\xymatrix{ \Omega\Sigma X^{(4)}\ar[d]^{\Omega \Phi}\ar[rr]^{\Omega
W_4} & &\Omega\Sigma X\ar[d]^{H_2}\\
\Omega\Sigma X^{(4)}\ar[rr]^{\Omega W_2(X^{(2)})} && \Omega\Sigma
X^{(2)}.}
\]
\end{thm}
\begin{proof} The second Hopf invariant of the $4$--fold Whitehead
product
\[
\ddr{\Omega\Sigma X^{(4)}}{\Omega W_4}{\Omega\Sigma X}
{H_2}{\Omega\Sigma X^{(2)}}
\]
restricted to $X^{(4)}$ has as its algebraic counterpart the
coalgebra map
\[
\ddr{V^{\otimes 4}}{\beta_4}{T(V)}{H_2}{T(V^{\otimes 2})}.
\]
Recall that
\[
\begin{array}{ll}
\beta_4(x_1x_2x_3x_4)= &
x_1x_2x_3x_4-x_2x_1x_3x_4-x_3x_1x_2x_4+x_3x_2x_1x_4\\
& -x_4x_1x_2x_3+x_4x_1x_2x_3+x_4x_3x_1x_2-x_4x_3x_2x_1.
\end{array}
\]
According to the Cohen-Taylor combinatorial formulae~\cite{CT} for
the Hopf invariant
\[
H_2(a_1a_2a_3a_4)=a_1a_2a_3a_4 +a_1a_3a_2a_4 +a_2a_3a_1a_4,
\]
the direct calculation gives
\[
H_2\circ\beta_4(x_1x_2x_3x_4)=[x_1x_2,x_4x_3]-[x_2x_1,x_4x_3]+[x_4x_1,x_3x_2]+[x_3x_1,x_4x_2].
\]
By taking the geometrical realization of the above formulae, we
obtain the map $S_2(X^{(2)})\circ\Phi$, where $S_2(X^{(2)})$ is the
$2$--fold Samelson product. Notice that the Samelson product is the
adjoint of the Whitehead product. Hence by taking the (unique)
multiplicative extension, the assertion follows.
\end{proof}
\section{Application: The rate of exponent growth for the Barratt Conjecture}
\label{sec4}

In this section we are concerned with an application of
combinatorial method to the exponent problem in homotopy theory.
Closely related to the exponent problem is the Barratt Conjecture
stated as follows. If $f\co \Sigma^2X\lra Y$ is of order $p^r$ in
$[\Sigma^2X,Y]$, then
\[
p^{r+1}\Im \left(f_*\co\pi_*(\Sigma^2 X)\lra\pi_*(Y)\right)=0.
\]
In particular, if the identity map on $\Sigma^2 X$ has order $p^r$
in $[\Sigma^2X,\Sigma^2X]$, then the exponent of $\Sigma^2X$ is
$p^{r+1}$, that is,
\[
p^{r+1}\pi_*(\Sigma^2X)=0.
\]
A stronger version of the Barratt Conjecture is concerned with the
multiplicative exponent of $\Sigma^2X$ and can be stated as follows.
If $f\co \Sigma^2X\lra Y$ is of order $p^r$ in $[\Sigma^2X,Y]$,
then
\[
\Omega^2f\co \Omega^2\Sigma^2X\lra\Omega^2Y
\]
has an order bounded by $p^{r+1}$ in
$[\Omega^2\Sigma^2X,\Omega^2Y]$.

An early result due to M\,G Barratt~\cite{Barratt} gave a bound on
the rate of the exponent growth for the Barratt Conjecture. The next
theorem states the analogue rate of the exponent growth but for the
stronger version of the Barratt Conjecture. In this section we
assume that $R=\Z/p^r$.
\begin{thm}\label{prop3.12.1}
Let $X=\Sigma X'$ be a suspension and let $f\co X\lra \Omega Y$
be a map such that $p^r[f]=0$ in the group $[X,\Omega Y]$. Let
$J(f)\co J(X)\lra \Omega Y$ be the canonical multiplicative
extension of $f$. Then the following hold.
\begin{enumerate}
\item The map
\[
J(f)|_{J_n(X)}\co J_n(X)\lra \Omega Y
\]
has order $p^{r+t}$ in $[J_n(X),\Omega Y]$ if $n<p^{t+1}$.
 \item The composite
\[
\ddr{J_{p^{t+1}}(X)}{J(f)|_{J_{p^{t+1}}(X)}}{\Omega
Y}{p^{r+t}}{\Omega Y}
\]
is homotopic to the composite
\[
\ddddr{J_{p^{t+1}}(X)}{\pinch}{X^{(p^{t+1})}}
{p^{r-1}(\sum_{\tau\in\Sigma_{p^{t+1}-1}}1\wedge
\tau)}{X^{(p^{t+1})}}{\widetilde{W}_{p^{t+1}}}{J(X)}{J(f)}{\Omega
Y,}
\]
where $p^{r+t}\co \Omega Y\lra \Omega Y$ is the $p^{r+t}$-th
power map, $\widetilde{W}_n$~is the $n$--fold Samelson product and
\[
1\wedge \tau (x_1\wedge \cdots\wedge x_{p^{t+1}})=x_1\wedge
x_{\tau(2)}\wedge\cdots\wedge x_{\tau(p^{t+1})}\co
X^{(p^{t+1})}\lra X^{(p^{t+1})}
\]
is the map which permutes positions. \item Let $g=J(f)\circ
\widetilde{W}_{p^{t+1}}\circ
(\sum_{\tau\in\Sigma_{p^{t+1}-1}}1\wedge\tau)\circ p^{r-1}\co
X^{(p^{t+1})}\lra\Omega Y$. Then $g$ is an equivariant map with
respect to the symmetric group action, that is,
\[
g\circ\sigma\simeq g
\]
for any $\sigma\in \Sigma_{p^{t+1}}$.
\end{enumerate}
\end{thm}
\begin{rem}
The map
\[
p^{r-1}(\sum_{\tau\in\Sigma_{p^{t+1}-1}}1\wedge \tau)\co
X^{(p^{t+1})}\lra X^{(p^{t+1})}
\]
is well-defined because $X$~is a suspension and so
$[X^{(p^{t+1})},X^{(p^{t+1})}]$ is an abelian group.
\end{rem}
\begin{proof}
Recall that the ground ring $R$ is $\Z/p^r\Z$. As $f\co X\lra
\Omega Y$ has order $p^r$, \fullref{cor:H_n^qtoM} implies that
the homomorphism $e_f\co\H\lra[J(X), \Omega Y]$ factors through
the quotient group $\H^{\Z/q}$, where $q=p^r$. There is an induced
homomorphism $e^{\Z/q}_f\co\H^{\Z/q}\lra[J(X), \Omega Y]$.
Further, recall that there is a group isomorphism
$e\co\H^{\Z/q}\lra \Coalg^{\Z/q}(T(-),T(-))$.

To show that the exponent of $J_n(f)$ is $p^{r+t}$ for $n< p^{t+1}$,
consider first the counterpart of the $p^{r+t}$-th power map
$p^{r+t}\co J(X)\lra J(X)$ on the level of natural coalgebra
transformations of the tensor algebra $T(-)$, that is, the
$p^{r+t}$-th convolution power of the identity $\Id_T$
\[
\ddr{\phi_t= \Id ^{*p^{r+t}}\co T}{\Psi_{p^{r+t}-1}}{T^{\otimes
p^{r+t}}}{\mu}{T}.
\]
Note that $\Hom_{R}(T,T)$ is an algebra under the convolution
product with \linebreak $\Coalg( T,T)$ $\subseteq \Hom_R(T,T)$.
The identity $1$ in the ring $\Hom_R(T,T)$ is the composite
\[
\ddr{1\co T}{\epsilon}{R}{\nu}{T}.
\]
Let $\overline\Id=\Id_T-1$ in $\Hom_R(T,T)$ which is represented by
the composite
\[
\overline\Id\co T\lra IT\lra T.
\]
Then
\[
\Id^{*p^{r+t}}=(1+\overline{\Id})^{p^{r+t}}
=1+\sum_{k=1}^{p^{r+t}}\binom{p^{r+t}}{k} \overline{\Id}^{*k}
\]
\[
=1+\binom{p^{r+t}}{p^t}\overline{\Id}^{*p^t}+\sum_{k=p^t+1}^{p^{r+t}}\binom{p^{r+t}}{k}
\overline{\Id}^{*k}
\]
in $\Hom_R(T,T)$ because $\binom{p^{r+t}}{i}\equiv 0$ mod $p^r$ for
$1<i<p^t$. The $k$-th convolution power $\overline{\Id}^{*k}$ is
represented by the composite
\[
\dddr{T}{\Psi_{k-1}}{T^{\otimes k}}{}{IT^{\otimes k}}{\mu}{T}.
\]
Now the restriction map $\overline{\Id}^{*k}\mid _{J_{p^t}}=0$ for
$k>p^t$ because $IT^{\otimes k}$ is a summand of $T_k$ and a
collection of $T_j$'s with $j>k$. Moreover,
$\overline{\Id}^{*p^t}\mid _{J_{p^t}}$ is represented by the
composite
\[
\ddr{J_p(V)}{}{T_{p^t}(V)=V^{\otimes p^t}}{\sum_{\sigma\in
\Sigma_{p^t}}}{V^{\otimes p^t}}\lra T(V),
\]
where $\Sigma_{p^t}$ acts on $V^{\otimes p^t}$ by permuting
positions. It follows that
\[
\overline{\Id}_T^{*p^t}\in \ker \left(\Coalg(J_{p^t},T)\lra
\Coalg(J_{p^t-1},T)\right)\cong\Lie^R(p^t)
\]
represented by the element $\binom{p^{r+t}}{p^t}\sum_{\sigma\in
\Sigma_{p^t}}\sigma$. We need to rewrite this element in terms of
Lie elements. Let $\bar{V}_n$ be the free $R$--module with a basis
$\{ x_1,\ldots, x_n\}$ and let $\gamma_n^R$ be the $R$--submodule
of $\bar V^{\otimes n}$ generated by the homogenous elements
$x_{\sigma(1)}x_{\sigma(2)}\negthinspace\cdots x_{\sigma(n)}$ for
$\sigma\in S_n$. Let $\Lie^R(n)$ be the $R$--submodule of
$\gamma^R_n$ generated by the $n$--fold commutators
$[[x_{\sigma(1)}, x_{\sigma(2)},\cdots x_{\sigma(n)}]$ for
$\sigma\in S_n$. Let
\[
\tr_n=\sum_{\sigma\in \Sigma_n} x_{\sigma(1)}x_{\sigma(2)}\cdots
x_{\sigma(n)}\in \gamma_n\subseteq \bar{V}^{\otimes n} \quad
\text{and}
\]
\[
\overline{\tr}_n=\sum_{\tau\in
\Sigma_{n-1}}[[x_1,x_{\tau(2)}],x_{\tau(3)},\cdots,x_{\tau(n)}]\in\Lie^R(n)
\]
be the sum of the standard basis for $\gamma_n$ and $\Lie^R(n),$
respectively. As $R=\Z/p^r$,
\[
p^{r+t-1}(x_1+x_2\cdots +x_{p^l})^{p^l}=p^{r+t-1}\tr_{p^l} + W
\]
\[
=p^{r+t-1}(x_1+v)^{p^l}=p^{r+t-1}[[x_1,v],v,\ldots
v]+W'=p^{r+t-1}\overline{\tr}_{p^l}+W
\]
where $ v=x_2+x_3+\cdots x_{p^l}$, $W$ is a sum of the homogeneous
terms in which one of the $x_i$'s occurs at least twice and $W'$ is
a sum of the homogeneous terms in which the number of occurrences of
$x_1$ is $0$ or greater than $1$. Thus
\[
p^{r+t-1}\tr_{p^l}=p^{r+t-1} \overline{\tr}_{p^l}
\]
and so
\[
\Id_T^{*p^{r+t}}=\binom{p^{r+t}}{p^t}\overline{\tr}_{p^t}
\]
in $\Lie^R(p^r)$.
\medskip

Being an element of the Cohen group $\overline{\tr}_n$ has a
geometrical realization given by the composite
\[
\ddddr{\overline{\tr}\co J_n(\Sigma X)}{\rm{pinch}}{(\Sigma
X)^{(n)}}{\sum_{\tau\in \Sigma_{n-1}}\id_{\Sigma X}\wedge
\tau}{(\Sigma X)^{(n)}}{W_n}{\Omega\Sigma^2X}{\Omega f}{\Omega Y}
\]
for $f\co \Sigma^2X\lra Y$.

Now consider the case when $n=p^{t+1}$. Then
\[
\begin{array}{l} e((x_1\cdots x_{p^t})^{p^{r+t}})=1+ (\begin{array}{ll}
p^{r+t}\\
p^{t+1}\\
\end{array})
\sum_{\sigma\Sigma_{p^t}}y_{\sigma(1)}\cdots y_{\sigma(p^t)} \\
=1+p^{r-1}q\sum_{\sigma\in\Sigma_{p^t}}y_{\sigma(1)}\cdots
y_{\sigma(n)}.\\
\end{array}
\]
Cohen, Moore and Neisendorfer~\cite{CMN} proved that
\[
\sum_{\tau\in\Sigma_{p^t-1}}[[y_1,y_{\tau(2)},\cdots,y_{\tau(p^t)}]
\equiv\sum_{\sigma\in\Sigma_{p^t}}y_{\sigma(1)}\cdots
y_{\sigma(n)}\pmod{p},
\]
where $\Sigma_{p^t-1}$ acts on $\{2,\cdots,p^t\}$. Thus
\[
p^{r-1}\sum_{\tau\in\Sigma_{p^t-1}}[[y_1,y_{\tau(2)},\cdots,y_{\tau(p^t)}]
=p^{r-1}\sum_{\sigma\in\Sigma_{p^t}}y_{\sigma(1)}\cdots
y_{\sigma(n)}.
\]
Assertions (2) and (3) follow.
\end{proof}

For $X$ a $p$--local suspension space, we now consider the special map
$$\phi_n=\sum_{\sigma\in\Sigma_n}\sigma\co X^{(n)}\lra X^{(n)},$$
where $\Sigma_n$ acts on $X$ by permuting positions.
\begin{lem}
\label{lem:tracfact} Let $\alpha\in\Z_{(p)}(\Sigma_n)$, such that
$\alpha=\sum_{\sigma\in\Sigma_n}k_{\sigma}\sigma$. Then
\[
\phi_n\circ\alpha=\left(\sum_{\sigma\in\Sigma_n}k_{\sigma}\right)\phi_n.
\]
\end{lem}
\begin{proof}
The statement follows from the observation that
$\phi_n\circ\sigma=\phi_n$ for any permutation $\sigma\in\Sigma_n$.
\end{proof}
\begin{lem}\label{lem3.14}
Let $f=\sum_{\sigma\in\Sigma_n}k_{\sigma}\sigma\in\Z(\Sigma_n)$,
where $k_{\sigma}\in\Z$, such that the sum of the coefficients
\[
\chi(f)=\sum_{\sigma\in\Sigma_n}k_{\sigma}\not\equiv 0\pmod{p}.
\]
Let $X$ be a $p$--torsion suspension of finite type. Then there
exists a map
\[
\phi_n(f)\co\hocolim_{f}X^{(n)}\lra\hocolim_{f}X^{(n)}
\]
such that the diagram
\[
\xymatrix{ X^{(n)}\ar[r]^{\phi_n}\ar[d]^{\hocolim(f)}&X^{(n)}\\
\hocolim_fX^{(n)}\ar[r]^{\phi_n(f)}&\hocolim_fX^{(n)}\ar[u]^{s}}
\]
commutes up to homotopy, where $s\co\hocolim_f X^{(n)}\lra
X^{(n)}$ is a map such that $\hocolim(f)\circ
s\co\hocolim_fX^{(n)}\lra\hocolim_f X^{(n)}$ is the identity map.
\end{lem}
\begin{proof}
We may assume that $\chi(f)=1$. Notice that
\[
\phi_n=f^q\circ \phi_n\circ f^r
\]
for any $q, r\geq 1$, where $f^m$ is the $m$--fold composition
$f\circ f\circ\cdots\circ f$. The assertion follows.
\end{proof}
Assume that $X$ is a suspension. Let $M_n(X)$ be the smallest
functorial retract of $X^{(n)}$ that contains the bottom cell in the
sense of papers~\cite{SELICKWU,CSW}.

\begin{proposition}
\label{trace} Let $X$ be a suspension. Then the trace map
\[
\phi_n=\sum_{\sigma\in\Sigma_n}\sigma\co X^{(n)}\lra X^{(n)}
\]
factors through the smallest functorial retract $M_n(X)$ of
$X^{(n)}$.
\end{proposition}
\begin{proof}
By construction of $M_n(X)$, there exists an idempotent $e_n$ in the
group algebra $\Z_{(p)}(\Sigma_n)$ such that
$e_n=\sum_{\sigma\in\Sigma_n}k_{\sigma}\sigma$ with
$\sum_{\sigma\in\Sigma_n}k_{\sigma}=1$ and the induced geometric map
$e_n\co X^{(n)}\lra X^{(n)}$ has the homotopy colimit
\[
M_n(X)=\mathrm{hocolim}_{e_n}X^{(n)}.
\]
Now the result follows from \fullref{lem:tracfact}.
\end{proof}
\begin{corollary}\label{cor:factorization}
Let $X$ be a suspension. Then the composite
\[
\ddr{J_{p^{t+1}}(X)}{J(f)|_{J_{p^{t+1}}(X)}}{\Omega
Y}{p^{r+t}}{\Omega Y}
\]
is homotopic to the composite
\[
\ddr{J_{p^{t+1}}(X)}{\pinch}{X^{(p^{t+1})}}{}{X\wedge
M_{p^{t+1}-1}(X)}\lra\ddr{X^{(p^{t+1})}}{W_{p^{t+1}}}{J(X)}{J(f)}{\Omega
Y}.
\]
\end{corollary}

The homology of $M_n(X)$ is unknown for a general $X$. The
determination of the homology of $M_n(X)$ is equivalent to the
fundamental problem in the modular representation theory of the
symmetric group $\Sigma_n$, that of determining the decomposition of
the group ring into indecomposable summands.

In the case of two cell complexes the homology of $M(X)$ is well
understood. Let $S_n(V)$ denote the $n$-th homogeneous component of
the symmetric algebra on $V$.
\begin{thm}{\rm\cite[Corollary 1.5]{CSW}}\qua
\label{thm:trace2cell}Let $X$ be a two-cell suspension. Then
\[
\tilde{H}_*(M_n(X))\cong S_n(\tilde{H}_*(X))
\]
provided that $n=cp^r-1$ for some $1\leq c\leq p-1$.
\end{thm}

\fullref{thm:trace2cell} gives a special meaning to
\fullref{cor:factorization}, claiming that the first
obstruction to the Barratt Conjecture for a two cell complex $X$
factors through the very small space $X\wedge M_{p^{t+1}-1}(X)$.

Let us consider a particular example in which~$p=2$. Let $X$ be a
suspension. We write $[q]\co X\lra X$ for the co-H $q$-th power
map for any integer $q$.

\begin{proposition}\label{proposition3.17}
Let $X=P^n(2)$ be the $n$--dimensional mod $2$ Moore space with
$n\geq 3$. Then the composite
\[
\ddr{X^{(3)}}{[2]}{X^{(3)}}{\sum_{\sigma\in\Sigma_3}\sigma}{X^{(3)}}
\]
is null homotopic and therefore by the
\fullref{application_Barrattintro}~$(2)$, the power map
\[
8|_{J_4(P^n(2))}\co J_4(P^n(2))\lra J(P^n(2))\simeq \Omega
(P^{n+1}(2))
\] is null homotopic.
\end{proposition}

\begin{proof}
Let $(123)$ be the $3$--cycle in $\Sigma_3$ and let
$f=(123)+(123)^2+(123)^3$. Then one can check that
$\hocolim_fX^{(3)}\simeq \mathbb{C}P^2\wedge P^{3n-4}(2)$, as
in~\cite{Wu4}. By \fullref{lem3.14}, there exists a map
$\phi_3(f)\co \mathbb{C}P^2\wedge P^{3n-4}(2)\lra
\mathbb{C}P^2\wedge P^{3n-4}(2)$ such that the map
$\sum_{\sigma\in\Sigma_3}\sigma\co X^{(3)}\lra X^{(3)}$ is
homotopic to the composite
\[
\dddr{X^{(3)}}{}{\mathbb{C}P^2\wedge
P^{3n-4}(2)}{\phi_3(f)}{\mathbb{C}P^2\wedge P^{3n-4}(2)}{}{X^{(3)}.}
\]
Notice that
\[
\phi_3(f)_*\co H_{3n-3}(\mathbb{C}P^2\wedge P^{3n-4}(2);\Z/2)\lra
H_{3n-3}(\mathbb{C}P^2\wedge P^{3n-4}(2);\Z/2)
\]
is zero. Thus $\phi_3(f)$ restricted to the bottom cell is null
homotopic. Notice that $\pi_{3n-2}(\mathbb{C}P^2\wedge
P^{3n-4}(2))=0$. Thus $\phi_3(f)$ restricted to the
$(3n-2)$--skeleton of the $4$--cell complex $\mathbb{C}P^2\wedge
P^{3n-4}(2)$ is null homotopic. Thus there exists a map
$\bar\phi\co P^{3n}(2)\lra \mathbb{C}P^2\wedge P^{3n-4}(2)$ such
that $\phi_3(f)$ is homotopic to the composite
\[
\ddr{\mathbb{C}P^2\wedge
P^{3n-4}(2)}{q}{P^{3n}(2)}{\bar\phi}{\mathbb{C}P^2\wedge
P^{3n-4}(2),}
\]
where $q$ is the pinch map. Notice that the map $[2]\co
P^{3n}(2)\lra P^{3n}(2)$ is homotopic to the composite
\[
\dddr{P^{3n}(2)}{\pinch}{S^{3n}}{\eta}{S^{3n-1}}{}{P^{3n}(2).}
\]
Thus the pinch map $q\co \mathbb{C}P^2\wedge P^{3n-4}(2)\lra
P^{3n}(2)$, which is a suspension map, is of exponent $2$ in the
group $[\mathbb{C}P^2\wedge P^{3n-4}(2),P^{3n}(2)]$. The assertion
follows.
\end{proof}

\bibliographystyle{gtart}
\bibliography{link}

\begin{thebibliography}{}
\providecommand\bibmarginpar{\leavevmode\marginpar}
\def\urlstyle#1{{\tt #1}}

\bibitem{Barratt}
\textbf{M\,G Barratt}, \emph{Spaces of finite characteristic}, Quart. J. Math.
  Oxford Ser. $(2)$ 11 (1960) 124--136 \xox{MR}{0120647}

\bibitem{Cohen}
\textbf{F\,R Cohen}, \emph{On combinatorial group theory in homotopy theory},
  preprint

\bibitem{CMN}
\textbf{F\,R Cohen}, \textbf{J\,C Moore}, \textbf{J\,A Neisendorfer},
  \href{http://links.jstor.org/sici?sici=0003-486X(197901)2:109:1%3C121:TIHG%3%
E2.0.CO%3B2-J} {\emph{Torsion in homotopy groups}}, Ann. of Math. $(2)$ 109
  (1979) 121--168 \xox{MR}{519355}

\bibitem{CSW}
\textbf{F\,R Cohen}, \textbf{P Selick}, \textbf{J Wu}, \emph{Natural
  decompositions of self-smashes of 2-cell complexes}, preprint

\bibitem{CT}
\textbf{F\,R Cohen}, \textbf{L\,R Taylor},
  \href{http://dx.doi.org/10.1007/BF01184668} {\emph{Homology of function
  spaces}}, Math. Z. 198 (1988) 299--316 \xox{MR}{946606}

\bibitem{GW}
\textbf{J Grbi\'c}, \textbf{J Wu},
  \href{http://dx.doi.org/10.2140/agt.6.2006.2189} {\emph{Natural
  transformations of tensor algebras and representations of combinatorial
  groups}}, Algebr. Geom. Topol. 6 (2006) 2189--2228

\bibitem{N}
\textbf{J\,A Neisendorfer}, \emph{The exponent of a {M}oore space}, from:
  ``Algebraic topology and algebraic $K$-theory (Princeton, N.J., 1983)'', Ann.
  of Math. Stud. 113, Princeton Univ. Press (1987)  35--71 \xox{MR}{921472}

\bibitem{SW}
\textbf{P Selick}, \textbf{J Wu}, \emph{On natural coalgebra decompositions of
  tensor algebras and loop suspensions}, Mem. Amer. Math. Soc. 148 (2000)
  \xox{MR}{1706247}

\bibitem{SELICKWU}
\textbf{P Selick}, \textbf{J Wu},
  \href{http://dx.doi.org/10.1007/s00229-002-0353-1} {\emph{On functorial
  decompositions of self-smash products}}, Manuscripta Math. 111 (2003)
  435--457 \xox{MR}{2002820}

\bibitem{Th}
\textbf{S Theriault}, \emph{Exponents of mod-$2^r$ Moore spaces}, to appear in
  Topology

\bibitem{Wu4}
\textbf{J Wu}, \emph{On maps from loop suspensions to loop spaces and the
  shuffle relations on the Cohen groups}, to appear in Memoirs AMS

\bibitem{Wu}
\textbf{J Wu}, \href{http://dx.doi.org/10.1016/S0040-9383(97)00058-X} {\emph{On
  combinatorial calculations for the {J}ames-{H}opf maps}}, Topology 37 (1998)
  1011--1023 \xox{MR}{1650426}

\end{thebibliography}

\end{document}